\documentclass[12pt, reqno]{amsart}

\usepackage{graphicx}
\usepackage{hyperref}

\usepackage{vmargin}
\setpapersize{USletter}
\setmargrb{1in}{1in}{1in}{1in}

\newcommand{\KK}{\mathbb K} 

\newcommand{\LL}{\mathbb L} 

\newcommand{\R}{\mathbb R}
\newcommand{\Q}{\mathbb Q}
\newcommand{\Z}{\mathbb Z}
\newcommand{\N}{\mathbb N}
\newcommand{\PP}{\mathbb P} 
\newcommand{\Hirz}{\mathbb{F}_2} 

\renewcommand{\O}{\mathcal O} 

\newcommand{\bA}{\mathbf A}
\newcommand{\bB}{\mathbf B}
\newcommand{\bC}{\mathbf C}
\newcommand{\bD}{\mathbf D}
\newcommand{\bE}{\mathbf E}
\newcommand{\bF}{\mathbf F}
\newcommand{\bG}{\mathbf G}

\newcommand\eu{E_1}
\newcommand\ed{E_2}
\newcommand\et{E_3}
\newcommand\eq{E_4}
\newcommand\lud{L_{12}}
\newcommand\lut{L_{13}}
\newcommand\luq{L_{14}}
\newcommand\ldt{L_{23}}
\newcommand\ldq{L_{24}}
\newcommand\ltq{L_{34}}

\newcommand{\RYeff}{R^+_Y} 
\newcommand{\money}{\mathcal{C}} 

\newcommand{\defining}[1]{\textbf{#1}}

\newcommand{\inters}[2]{\langle{#1},{#2}\rangle}

\newcommand{\rt}{x}

\newcommand{\orbit}{\mathcal{O}}
\newcommand{\orbitrootsystem}[2]{\orbit(#1,#2)}
\newcommand{\orbsys}{\orbitrootsystem}

\newcommand{\dualcone}[1]{#1^{\vee}}
\newcommand{\GalK}{G_{\KK}}

\DeclareMathOperator{\Nef}{Nef}

\DeclareMathOperator{\Eff}{Eff}

\DeclareMathOperator{\Pic}{Pic}

\DeclareMathOperator{\Vol}{Vol}
\DeclareMathOperator{\Bl}{Bl}
\DeclareMathOperator{\Gal}{Gal}

\DeclareMathOperator{\id}{id} 
\DeclareMathOperator{\linspan}{span}
\DeclareMathOperator{\Aut}{Aut}

\theoremstyle{plain}
\newtheorem{thm}{Theorem}[section]
\newtheorem{lem}[thm]{Lemma}
\newtheorem{cor}[thm]{Corollary}
\newtheorem{prop}[thm]{Proposition}
\newtheorem*{thm*}{Theorem}

\theoremstyle{definition}
\newtheorem{defn}[thm]{Definition}

\newtheorem{convention}[thm]{Convention}
\newtheorem{example}[thm]{Example}

\theoremstyle{remark}
\newtheorem{rem}[thm]{Remark}

\newtheoremstyle{citing}
 {}
 {}
 {\itshape}
 {}
 {\bfseries}
 {.}
 {.5em }
 {\thmnote{#3}}

\theoremstyle{citing}
\newtheorem*{varthm}{} 

\title{The nef cone volume of generalized Del Pezzo surfaces}

\author{Ulrich Derenthal}
\address{Institut F\"{u}r Mathematik, Universit\"{a}t Z\"{u}rich,
Winterthurerstrasse 190, 8057 Z\"{u}rich, Switzerland}
\email{ulrich.derenthal@math.unizh.ch}

\author{Michael Joyce}
\address{Department of Mathematics, Tulane University, New Orleans LA 70118}
\email{mjoyce@math.tulane.edu}

\author{Zach Teitler}
\address{Department of Mathematics, SLU 10687, Hammond, LA 70402}
\email{zteitler@selu.edu}

\date{7/24/07}



\begin{document}

\bibliographystyle{amsalpha}

\begin{abstract}
  We compute a naturally defined measure of the size of the nef cone of a Del
  Pezzo surface.  The resulting number appears in a conjecture of Manin on the
  asymptotic behavior of the number of rational points of bounded height on
  the surface.  The nef cone volume of a Del Pezzo surface $Y$ with
  $(-2)$-curves defined over an algebraically closed field is equal to the nef
  cone volume of a smooth Del Pezzo surface of the same degree divided by the
  order of the Weyl group of a simply-laced root system associated to the
  configuration of $(-2)$-curves on $Y$.
  When $Y$ is defined over a non-closed field of characteristic $0$,
  a similar result holds, except that the associated root system is no
  longer necessarily simply-laced.
\end{abstract}

\maketitle

\section{Introduction}\label{intro}

An \defining{ordinary Del Pezzo surface} is a smooth projective rational
surface $X$ on which the anticanonical class $-K_X$ is ample. If $X$ is
defined over an algebraically closed field, then $X$ is one of the following:
$\PP^2$, $\PP^1 \times \PP^1$, or the blowup of $\PP^2$ at up to $8$
points in general position.  Points are in general position if no three are
collinear, no six lie on a conic, and no eight lie on a cubic with one of them
a singular point of the cubic. Then $X$ may contain $(-1)$-curves, but no
$(-2)$-curves, where for $n \in \{1,2\}$, a \defining{$(-n)$-curve} is
a smooth rational curve on $X$ which has self intersection number $-n$.

A \defining{generalized Del Pezzo surface} is a smooth projective rational
surface $Y$ on which $-K_Y$ is big and nef.
If $Y$ is defined over an algebraically closed field,
then $Y$ is one of the following: $\PP^2$, $\PP^1 \times \PP^1$,
the Hirzebruch surface $\Hirz$,
or a surface obtained from $\PP^2$ by a sequence
of blowings-up at up to $8$ points, possibly infinitely near,
each not lying on any $(-2)$-curve.
Over an algebraically closed field,
a generalized Del Pezzo surface is ordinary if and only if it contains
no $(-2)$-curves.  See Section~\ref{section:gdp} for more details.

The nef cone volume of a generalized Del Pezzo surface $Y$ is
equal to the volume of the cross-section of the nef cone of $Y$
obtained by intersecting with the hyperplane consisting of those
divisor classes whose intersection with the anticanonical class
$-K_Y$ is equal to $1$. The resulting cross-section is a polytope.
Its volume is a rational number, denoted $\alpha(Y)$. We give more
details of this definition in Section~\ref{section:nef cone volume}.

In this paper we compute $\alpha(Y)$ for generalized Del Pezzo surfaces $Y$.

The \defining{degree} of $Y$ is the self-intersection number
$d=\inters{-K_Y}{-K_Y}$.  The degree is always in the range $1 \leq d \leq 9$,
and when $Y$ is the blowup of $\PP^2$ at $r$ points in almost general
position, $d = 9 - r$.

A generalized Del Pezzo surface $Y$ defined over a field $\KK$ is
\defining{split} if it is either $\PP^2$, $\PP^1 \times \PP^1$,
$\mathbb{F}_2$, or the blow-up of $\PP^2$ at $1 \leq r \leq 8$ $\KK$-rational
points in almost general position.  See Definition~\ref{defn:split}.
Otherwise, $Y$ is said to be non-split (e.g., the blow-up of $\PP^2$ at two
conjugate points).  We consider only split $Y$ until Section
\ref{section:galois}, and so the reader may assume that $\KK$ is algebraically
closed until that point.

An investigation of $\alpha(Y)$ for split ordinary Del Pezzo surfaces
was undertaken by the first author.
He proved the following result:

\begin{thm}\cite[Theorem~4]{derenthal-alpha}\label{thm:derenthal}
  Let $X_d$ denote a split ordinary Del Pezzo surface of degree $d$ obtained
  by blowing up $9 - d$  points in general position on $\PP^2$ and let
  $N_d$ denote the number of $(-1)$-curves on $X_d$.  For $d \leq 6$,
  \begin{equation*}
    \alpha(X_d) = \frac{N_d}{d(9-d)} \alpha(X_{d+1}).
  \end{equation*}
\end{thm}

Combining this with simple calculations that show $\alpha(\PP^2) =
1/3$, $\alpha(\PP^1 \times \PP^1) = 1/4$, $\alpha(X_8) = 1/6$, and
$\alpha(X_7) = 1/24$,
this theorem allows for an inductive calculation of
$\alpha(X)$ for any split ordinary Del Pezzo surface $X$. This calculation
is summarized in Table~\ref{table:alpha of ordinary DPs}.

\begin{table}\label{table:alpha of ordinary DPs}
\[\begin{array}{| c | c | c | c | c | c | c | c | c | c |  }
\hline
d & 8 & 7 & 6 & 5 & 4 & 3 & 2 & 1 \\
\hline
N_d & 1 & 3 & 6 & 10 & 16 & 27 & 56 & 240 \\
\hline
\alpha(X_d) & 1/6 & 1/24 & 1/72 & 1/144 & 1/180 & 1/120 & 1/30 & 1 \\
\hline
\end{array}\]
\caption{Values of $\alpha(X_d)$ for ordinary Del Pezzo surfaces
$X_d$}\label{table:alpha-ordinary-dps}
\end{table}

We extend this result in two directions.  First, we study split generalized
Del Pezzo surfaces.  In Section~\ref{section:inductive}, we prove the
following theorem by analyzing the nef cone of such a surface $Y$. It allows
us to compute $\alpha(Y)$ by induction on the rank of the N{\'e}ron-Severi
group of $Y$.

\begin{thm}\label{thm:inductive}
  Let $Y$ be a split generalized Del Pezzo surface of degree $d \leq 7$.  For
  each $E$ in the set $\money$ of $(-1)$-curves on $Y$, let $Y_E$ denote the
  split generalized Del Pezzo surface of degree $d + 1$ obtained by
  contracting $E$.  Then
  \begin{equation*}
    \alpha(Y) = \sum_{E \in \money} \frac{1}{d(9-d)} \alpha(Y_E).
  \end{equation*}
\end{thm}

As with Theorem~\ref{thm:derenthal}, using the additional
calculation that $\alpha(\Hirz) = 1/8$,
this theorem allows us to compute $\alpha(Y)$ for any split
generalized Del Pezzo surface $Y$.

The first author computed $\alpha(Y)$ for split generalized Del Pezzo surfaces
$Y$ of degree $d \geq 3$ directly, using computer programs to find
a triangulation of the nef cone case by case~\cite[Section~3]{derenthal-alpha}.
This numerical data led us to formulate the following theorem; see
Section~\ref{section:weyl-groups-volume} for its proof.

\begin{thm}\label{thm:weyl}
Let $Y$ be a split generalized Del Pezzo surface of degree $d \leq 7$
and let $X$ be a split ordinary Del Pezzo surface of the same degree.
Then
\begin{equation*}
\alpha(Y) = \frac{1}{\#W(R_Y)} \alpha(X),
\end{equation*}
where $W(R_Y)$ is the Weyl group of the root system $R_Y$
whose simple roots are the $(-2)$-curves on $Y$.
\end{thm}

When combined with Theorem~\ref{thm:derenthal} the computation of $\alpha(Y)$
for an arbitrary split generalized Del Pezzo surface $Y$ of any degree is
reduced to a determination of the $(-2)$-curves on the surface.  See
Section~\ref{section:weyl} for more information on the root system $R_Y$ and
its Weyl group.

We also consider the case of non-split surfaces.  Suppose that $Y$ is a
generalized Del Pezzo surface and $X$ is an ordinary Del Pezzo surface, both
of the same degree and defined over the same field $\KK$ of characteristic $0$.
Then the N{\'e}ron-Severi groups of $X$ and $Y$ coincide
(Proposition~\ref{prop:galois descent on neron-severi group}) and the absolute
Galois group of $\KK$ acts as a finite group $G$ of automorphisms of this
group (Proposition~\ref{prop:automorphisms of N^1}).  Assume that the Galois
actions associated to $X$ and $Y$ coincide.  The Galois action on the root
system $R_Y$ allows us to associate to $Y$ an orbit root system
$\orbsys{R_Y}{G}$ (Definition~\ref{dfn:orbit root system}).  Our third main
result is that under these assumptions
\begin{equation*}
\alpha(Y) = \frac{\alpha(X)}{\# W(\orbsys{R_Y}{G})}.
\end{equation*}
See Corollary~\ref{cor:gen-dp-weyl}.  The integer appearing in the denominator
is the order of a Weyl group and is straightforward to compute.  Thus all that
is left is to compute $\alpha(X)$.  There are a finite number of cases in each
degree $d$, one for each conjugacy class of subgroups of the Weyl group of a
canonically defined root system $R_d$ (Section~\ref{section:roots-on-dps}).
We perform the computations for $d \ge 5$ in Section~\ref{subsection:high
  degree}.

\subsection*{Manin's conjecture}

The primary motivation for our study of the nef cone volume $\alpha$ is its
appearance in Manin's conjecture on the number of rational points of bounded
height on Fano varieties defined over number fields, as described below.
Although the conjecture is now known not to hold for all Fano
varieties~\cite[Theorems~3.1--3.3]{BT2}, it has been verified in a large
number of cases, including some varieties for which the anticanonical class is
big but not ample.

Let $X$ be a smooth projective variety defined over a
number field $\KK$ for which $-K_X$ is big and
assume that the set $X(\KK)$ of rational points is Zariski dense.
Equip $X(\KK)$ with an anticanonical height function $H$ (consult
\cite[Part~B]{MS} for information on height functions) and for any
constructible set $U \subset X$ let
\begin{equation*}
\mathcal{N}_U(B) := \# \{ P \in U(\KK) : H(P) \leq B \}.
\end{equation*}
The original formulation of the Manin conjecture~\cite[Conjecture~B]{BM}
posits the existence of a Zariski open set $U \subset X$ such that for any
open set $V \subset U$
\begin{equation*}
  \mathcal{N}_V(B) \sim c(X) B (\log B)^{\rho - 1}
  \quad\text{asymptotically as $B \to \infty$,}
\end{equation*}
where $\rho$ is the N{\'e}ron-Severi rank of $X$. The conjecture
was initially made for Fano varieties, but a more ambitious
version of the conjecture relaxes the condition on $-K_X$ to
merely being big. The leading constant was given a conjectural
interpretation by Peyre~\cite[Definition~2.4]{Peyre} and Batyrev
and Tschinkel \cite[Theorem~4.4.4]{BT}.
They predict that
\begin{equation*}
c(X) = \alpha (X) \beta(X) \tau(X),
\end{equation*}
where $\alpha(X) \in \Q$ is the constant of interest in this
paper, $\beta(X) \in \N$ is a cohomological invariant of the
Galois action on the N{\'e}ron-Severi group of $X$, and $\tau(X)
\in \R$ is a volume of adelic points on $X$.

\subsection*{Acknowledgments}\label{acknowledge}

We thank Brian Harbourne and Brendan Hassett for helpful comments
and encouragement.

\section{Definition of the nef cone volume}\label{section:nef cone volume}

Here we recall the definition of $\alpha(X)$, first introduced by Peyre
in~\cite[Definition~2.4]{Peyre}.

Let $X$ be a smooth complete variety for which $-K_X$ is big.
We denote the intersection form on $X$ by $\inters \cdot \cdot$.
Recall that a divisor class $D$ on $X$ is \defining{numerically
trivial} if $\inters{D}{C}=0$ for all curves (equivalently, all
$1$-cycles) $C$ on $X$, and two divisor classes are
\defining{numerically equivalent} if their difference is
numerically trivial. One similarly defines numerical equivalence
of curves.
Numerical equivalence classes of divisors on $X$ form a
finitely-generated torsion-free abelian group $N^1(X)$ whose dual
group $N_1(X)$ consists of numerical equivalence classes of
1-cycles on $X$. Let $N^1(X)_{\R} = N^1(X) \otimes_{\Z} \R$ and
$N_1(X)_{\R} = N_1(X) \otimes_{\Z} \R$ be the associated Euclidean
spaces.
Inside $N^1(X)_{\R}$ lies the \defining{effective cone} $\Eff^1(X)$,
the closed convex cone spanned by the classes of effective divisors.

Recall that for a finite-dimensional real inner product space $V$ and a convex
cone $\Gamma \subset V$, the dual convex cone $\dualcone{\Gamma} \subset V$ is
defined by
\[ \dualcone{\Gamma} = \{ \, v \in V : \inters{v}{c} \geq 0 \text{ for all $c
  \in \Gamma$} \, \} . \] The cone $\Gamma^{\vee}$ is closed as a subspace of
the Euclidean space $V$. The dual $\dualcone{\Eff^1(X)}$ of the effective cone
of $X$ in $N^1(X)_\R$ is the movable cone of $X$
(see~\cite[Theorem~2.2]{math.AG/0405285}, \cite[Section~11.4.C]{pag2}).  Note
that when $X$ is a surface, $N_1(X) = N^1(X)$ and $\dualcone{\Eff^1(X)}$ is
the \defining{nef cone} of $X$, denoted $\Nef(X)$.

Since the cone $\dualcone{\Eff^1(X)}$ has infinite volume in $N_1(X)_{\R}$, a
natural means of measuring its ``size'' is to truncate the cone in an
(anti)canonical manner.  To do this, consider the hyperplane
\begin{equation*}
\mathcal{H}_X := \{ C \in N_1(X)_{\R} : \inters{-K_X}{C}  = 1 \} .
\end{equation*}
Note that since $-K_X$ is big by hypothesis, $\mathcal{H}_X$
intersects each ray of $\dualcone{\Eff^{1}(X)}$. We endow
$N_1(X)_{\R}$ with Lebesgue measure $ds$ normalized so that
$N_1(X)$ has covolume 1, and we endow $\mathcal{H}_X$ with the
induced Leray measure $d\mu$ with respect to the linear form
$\inters{-K_X}{\cdot}$. That is, letting $l$ be the linear form
$l(v) = \inters{-K_X}{v}$, we have $ds = d\mu \wedge dl$. We
construct the polytope
\begin{equation*}
\mathcal{P}_X := \dualcone{\Eff^1(X)} \cap \mathcal{H}_X
\end{equation*}
and define
\begin{equation*}
\alpha(X) := \Vol(\mathcal{P}_X) = \int_{\mathcal{P}_X} \, d\mu .
\end{equation*}

There are variants of this definition which differ only by a
dimensional factor. Let $\rho = \dim N_1(X)_{\R}$ and
\begin{equation*}
  \mathcal{C}_X := \{ C \in \dualcone{\Eff^1(X)} : \inters{-K_X}{C}
  \leq 1 \}
\end{equation*}
be the convex hull of $\mathcal{P}_X$ and the origin.  Then a
simple slicing argument shows that
\begin{equation*}
\alpha(X) = \rho \cdot \Vol(\mathcal{C}_X).
\end{equation*}  Additionally,
\begin{equation*}
\alpha(X) = \frac{1}{(\rho - 1)!}
\idotsint_{\dualcone{\Eff^1(X)}} \exp\left(-\inters{-K_X}{s}\right) ds,
\end{equation*}
with the bigness of $-K_X$ insuring the convergence of the
integral.

\begin{example}
Let us compute $\alpha(\PP^2)$.
We have $N^1(\PP^2)_{\R} \cong \R^1$,
with the real number $x \in \R$ corresponding to the (real) divisor class $xL$,
where $L$ is the class of a line in $\PP^2$.
Then the nef cone $\Nef(\PP^2) = \{ x \in \R : x \geq 0\}$
and the anticanonical class corresponds the real number $3$.
The hyperplane $\mathcal{H}_{\PP^2}$ is just $\{ 1/3 \}$.
The polytope $\mathcal{P}_{\PP^2}$ is also $\{ 1/3 \}$
and the convex hull $\mathcal{C}_{\PP^2} = [0,1/3]$.
Thus $\mathcal{C}_{\PP^2}$ has volume $1/3$
and so $\alpha(\PP^2) = 1 \cdot \Vol(\mathcal{C}_{\PP^2}) = 1/3$.
\end{example}

\begin{example}
As a second example, let $X_8$ be the blowup of $\PP^2$ at a
single point. Let $L$ be the class of the pullback of a line to
$X_8$ and let $E$ be the class of the exceptional divisor. Then
$N^1(X_8)$ is generated by $L$ and $E$. In $N^1(X_8)_{\R} \cong
\R^2$, with $(a,b)$ corresponding to $aL+bE$, the nef cone
$\Nef(X_8)$ is equal to $\{ (a,b) : a \geq 0, a+b \geq 0 \}$, that
is, the cone with extremal rays spanned by $L$ and $L-E$. The
anticanonical class corresponds to the point $(3,-1)$. The
hyperplane $\mathcal{H}_{X_8}$ is the line $3a + b = 1$. One
checks that $\mathcal{P}_{X_8}$ is the segment joining the points
$(1/3,0)$ and $(1/2,-1/2)$. Then $\mathcal{C}_{X_8}$ is the
triangle with vertices the above two points together with the
origin. The area of this triangle is $1/12$, and so $\alpha(X_8) =
2 \cdot \Vol(\mathcal{C}_{X_8}) = 1/6$.
\end{example}

\subsection*{Terminology}
Peyre~\cite{Peyre} introduced the notation $\alpha(X)$, but did not
give a name to this quantity.
We will refer to $\alpha(X)$ as the ``nef cone volume of $X$'' whenever
$X$ is a surface.

\section{Generalized Del Pezzo surfaces}\label{section:gdp}

As stated in the introduction, a \defining{generalized Del Pezzo surface} is a
smooth projective rational surface $Y$ on which $-K_Y$ is big and nef. If $Y$
is defined over an algebraically closed field, $Y$ is one of $\PP^2$, $\PP^1
\times \PP^1$, the Hirzebruch surface $\Hirz$, or $\PP^2$ blown up at $1 \leq
r \leq 8$ points in almost general
position~\cite[Definition~III.2.1]{Demazure}. To blow up $r$ points on $\PP^2$
in almost general position is to construct a sequence of morphisms
\begin{equation*}
Y = Y_r \rightarrow Y_{r-1} \rightarrow \cdots \rightarrow Y_1
\rightarrow Y_0 = \PP^2,
\end{equation*}
where each map $Y_i \rightarrow Y_{i-1}$ is the blow-up of
$Y_{i-1}$ at a point $p_i \in Y_{i-1}$ not lying on any
irreducible curves of
self-intersection number $-2$ on $Y_i$.

For $n \in \{1, 2\}$, a \defining{$(-n)$-class} on $Y$ is a divisor class $D$
such that $\inters D D = -n$ and $\inters{D}{-K_Y} = 2 - n$.  If such a class
is effective, then there is necessarily a unique curve in that class.  If this
curve is irreducible, we use the term \defining{$(-n)$-curve} both for this
curve and its class.  It follows from the genus formula that a $(-n)$-curve is
a smooth rational curve.  A simple calculation \cite[Tables~2,~3]{Demazure}
shows that the sets of $(-1)$- and $(-2)$-classes on a generalized Del Pezzo
surface are finite.

Let $Y$ be a generalized Del Pezzo surface defined over a field $\KK$. We
denote $Y \times_{\KK} \overline{\KK}$ by $\overline{Y}$.
Recall that a generalized Del Pezzo surface $Y$ is an ordinary Del Pezzo
surface if and only if the anticanonical class $-K_Y$ is ample.  Equivalently,
there are no $(-2)$-curves on $\overline{Y}$.

\begin{convention}
Throughout the paper, we will use $X$ to refer to an ordinary Del Pezzo surface
and $Y$ to refer to a generalized (possibly ordinary) Del Pezzo surface.
\end{convention}

The absolute Galois group $\GalK = \Gal(\overline{\KK}/\KK)$ acts
on $N^1(\overline{Y})$.

\begin{defn}\label{defn:split}
A generalized Del Pezzo surface $Y$ is \defining{split} if $Y(\KK)
\neq \emptyset$ and the action of $\GalK$ on $N^1(\overline{Y})$
is trivial.
\end{defn}

Excluding the exceptional cases where $\overline{Y}$ is isomorphic
to $\PP^1 \times \PP^1$ or $\mathbb{F}_2$, the existence of a
rational point assures that $Y$ is a blow-up of $\PP^2$ and the
triviality of the Galois action assures that each exceptional
divisor is defined over $\KK$, and thus the sequence of blown-up
points must themselves be defined over $\KK$.

In the remainder of this section, we prove that for a split $Y$,
the effective cone $\Eff^1(Y)$ is generated by the set of $(-1)$-
and $(-2)$-curves on $Y$, collecting a number of useful facts
along the way.

By the following lemma, the group $N^1(Y)$ depends only on the
degree of $Y$. We will make frequent use of this well-known
result.

\begin{lem}\label{lem:identification}
Let $X$ be a split ordinary Del Pezzo surface and let $Y$ be a
split generalized Del Pezzo surface of the same degree $d \leq 7$.
There is an isomorphism of $N^1(X)$ and $N^1(Y)$ which
identifies the intersection forms and takes $-K_X$ to $-K_Y$.
\end{lem}

\begin{proof}
Say $X$ is the blowup of $\PP^2$ at points $p_1,\dots,p_r \in \PP^2$, $r=9-d$,
with blowdown $\pi_X:X \to \PP^2$,
and say $Y$ is obtained by blowing up $\PP^2$ at points $q_1,\dots,q_r$:
\[ \pi_Y: Y=Y_r \to Y_{r-1} \to \cdots \to Y_1 \to Y_0=\PP^2 \] where $Y_{j} =
\Bl_{q_j}(Y_{j-1})$, $q_j \in Y_{j-1}$.  Let $E_{X,j}$ be the exceptional
divisor over $p_j$, and let $E_{Y,j}$ be the total transform in $Y$ of the
exceptional divisor over $q_j$.  (That is, if $f_j: Y \to Y_{j-1}$, then
$E_{Y,j} = f_j^{-1}(q_j)$, scheme-theoretically.)

Then $N^1(X)$ is the free abelian group on
$L_X=\pi_X^{*} \O_{\PP^2}(1)$, $E_{X,1}$, \dots, $E_{X,r}$.
Similarly, $N^1(Y)$ is the free abelian group on
$L_Y=\pi_Y^{*} \O_{\PP^2}(1)$, $E_{Y,1}$, \dots, $E_{Y,r}$.
The intersection form on $N^1(X)$ is given in this basis by
the diagonal matrix with entries $(1,-1,\dots,-1)$;
the intersection form on $N^1(Y)$ is given in this basis by the same matrix.
We have $-K_X = 3 L_X - \sum E_{X,j}$ and $-K_Y = 3 L_Y - \sum E_{Y,j}$.
\end{proof}

\begin{rem}
Note that the identification made in the proof of
Lemma~\ref{lem:identification} is not necessarily unique;
see~\cite[Theorem~0.1]{MR791295}.
\end{rem}

The next lemma is a modest generalization of \cite[Proposition~4.5]{HT}.

\begin{lem}\label{lem:gen cond}
Let $S$ be a surface and let $D_1, \dots, D_k$ be irreducible
effective divisors on $S$.  Let $\Gamma$ denote the cone generated
by $D_1, \dots, D_k$.  Then the effective cone of $S$ is equal to
$\Gamma$ if and only if $\dualcone{\Gamma} \subset \Gamma$.
\end{lem}

\begin{proof}
  If the effective cone of $S$ is equal to $\Gamma$ then it is a closed cone.
  The nef cone $\Nef(S) = \dualcone{\Gamma}$ is contained in the closure of
  the effective cone, which is just $\Gamma$.

  For the converse, it is clear that $\Gamma$ is contained in the effective
  cone of $S$.  Let $D$ be an effective divisor.  Then we can write $D = D' +
  a_1 D_1 + \cdots + a_k D_k$ with $a_i \geq 0$ and $D'$ having none of the
  $D_i$ as an irreducible component.  It is clear that $D'$ is contained in
  $\dualcone{\Gamma}$, and by hypothesis, $D'$ is consequently contained in
  $\Gamma$.  Hence the same is true of $D$.
\end{proof}

\begin{prop}\label{prop:exc are eff}
  If $Y$ is a split generalized Del Pezzo surface, every $(-1)$-class in
  $N^1(Y)$ is effective.  Indeed, if $E$ is any $(-1)$-class, then either
\begin{enumerate}
\item $E$ is a $(-1)$-curve;
\item $E$ can be written as the sum of a $(-1)$-curve and
one or more
  $(-2)$-curves;
\item $d = 1$ and $E$ can be written as the sum of $-K_Y$ and
one or more
  $(-2)$-curves.
\end{enumerate}
\end{prop}

\begin{proof}
See \cite[Theorem~III.2.c]{Demazure}.
\end{proof}

For a split generalized Del Pezzo surface $Y$ of degree $d \ge 2$, this shows
that every $(-1)$-class is a non-negative integral linear combination of
$(-1)$- and $(-2)$-curves. By the following lemma, this holds also in degree
$d=1$ if we allow rational instead of integral coefficients.

\begin{lem}\label{lem:degree 1}
  For a split generalized Del Pezzo surface $Y$ of degree~$1$, the
  anticanonical class $-K_Y$ is a linear combination of $(-1)$- and
  $(-2)$-curves with non-negative rational coefficients.
\end{lem}

\begin{proof}
  Let $X$ be an ordinary Del Pezzo surface of degree 1. It is easy to check
  that the sum of all $(-1)$-classes on $X$ is $-240K_X$.

  Using the identification of Lemma~\ref{lem:identification},
  the sum of all $(-1)$-classes on $Y$ is $-240K_Y$.
  Using Proposition~\ref{prop:exc are eff}, we can write $n$ of the
  $(-1)$-classes as the sum of a $(-1)$-curve and possibly some $(-2)$-curves,
  and the remaining $240-n$ of the $(-1)$-classes as the sum of $-K_Y$ and
  some $(-2)$-classes.
  Note that $E_{Y,8}$ in the proof of Lemma~\ref{lem:identification} is a
  $(-1)$-curve on $Y$, so we have $n>0$.

  This gives us $-240K_Y$ as the sum of $n$ $(-1)$-curves, $-(240-n)K_Y$, and
  some $(-2)$-curves. We transform this equation to write $-nK_Y$ as a sum of
  $(-1)$- and $(-2)$-curves.
\end{proof}

\begin{lem}\label{lem:irred of exc classes}
Let $Y$ be a split generalized Del Pezzo surface and let $E$ be a
$(-1)$-class in $N^1(Y)$.  Then $E$ is irreducible if and only if
$\langle E, C \rangle \geq 0$ for every $(-2)$-curve $C$.
\end{lem}

\begin{proof}
See \cite[Corollary on page 46]{Demazure}.
\end{proof}

In the case of ordinary Del Pezzo surfaces, the following result is well-known.

\begin{prop}\label{prop:effective cone generators smooth case}
Let $X$ be a split ordinary Del Pezzo surface of degree $d \leq
7$. Then the effective cone of $X$ is minimally generated by the
$(-1)$-classes on $X$, all of which are
$(-1)$-curves.
\end{prop}

\begin{proof}
  This can be proved directly (see~\cite[Theorem~V.4.11]{hartshorne} for a
  proof when $d=3$) or can be taken as an immediate consequence of the
  calculation of generators for the Cox ring given in~\cite[Theorem~3.2]{BP},
  making use of Lemma~\ref{lem:degree 1} in the case $d=1$.
\end{proof}

We now reach our main goal for this section.

\begin{thm}\label{thm:effective cone generators}
  If $Y$ is a split generalized Del Pezzo surface and has degree $d \leq 7$,
  the effective cone of $Y$ is finitely generated by the set of $(-1)$- and
  $(-2)$-curves.
\end{thm}

\begin{proof}
  Let $\Gamma$ be the cone generated by the $(-1)$- and $(-2)$-curves on $Y$.
  To prove the theorem, it suffices by Lemma~\ref{lem:gen cond} to show that
  $\dualcone{\Gamma} \subset \Gamma$. Let $X$ be a split ordinary Del Pezzo
  surface of the same degree as $Y$. Identify $N^1(X)$ and $N^1(Y)$ as in
  Lemma~\ref{lem:identification}.  Note that this identification takes
  $(-1)$-classes to $(-1)$-classes.  By Proposition~\ref{prop:effective cone
    generators smooth case}, $\Eff^1(X)$ is generated by $(-1)$-classes.  Each
  $(-1)$-class lies in $\Gamma$ by Proposition~\ref{prop:exc are eff} and
  Lemma~\ref{lem:degree 1}.  Therefore $\Eff^1(X) \subset \Gamma$.  It follows
  immediately that $\dualcone{\Gamma} \subset \dualcone{\Eff^1(X)}$.  From
  Lemma~\ref{lem:gen cond} we have $\dualcone{\Eff^1(X)} \subset \Eff^1(X)$.
  Thus $\dualcone{\Gamma} \subset \Gamma$ and hence $\Gamma = \Eff^1(Y)$,
  again by Lemma~\ref{lem:gen cond}.
\end{proof}

\begin{rem}\label{rem:lahayne-harbourne paper}
A generalization of Theorem 3.9 has already been proved by Lahyane
and Harbourne \cite[Lemma~4.1]{HL}.  We include our presentation
both as a summary of results that we will use later and also
because the approach here seems to have interest in its own right.
\end{rem}

\begin{cor}\label{cor:cone containments}
  Let $X$ be a split ordinary Del Pezzo surface and $Y$ a split generalized
  Del Pezzo surface with $\deg(X) = \deg(Y) \leq 7$.  Identifying $N^1(X)$ and
  $N^1(Y)$ as in Lemma~\ref{lem:identification}, we have $\Eff^1(X) \subset
  \Eff^1(Y)$ and $\Nef(X) \supset \Nef(Y)$.

  Let $\Gamma \subset N^1(Y)_{\R}$ be the cone spanned by the set of
  $(-2)$-curves on $Y$.  Then $\Eff^1(Y)$ is the sum of $\Eff^1(X)$
  and $\Gamma$, and $\Nef(Y) = \Nef(X) \cap \dualcone{\Gamma}$.
\qed
\end{cor}

\section{Inductive method}\label{section:inductive}

With these preliminaries in place, we now turn to proving
Theorem~\ref{thm:inductive}.
For a generalized Del Pezzo surface $Y$ and
any class $D \in N^1(Y)_{\R}$, we
denote by $D^{\perp}$ the hyperplane
\[ D^{\perp} := \{ C \in N_1(Y)_{\R} : \inters{D}{C} = 0 \} . \]

\begin{lem}\label{lem:nef cone identification}
Let $Y$ be a
split
generalized Del Pezzo surface and $E$ a $(-1)$-curve on $Y$.
Let $\pi: Y \to Y_E$ be the contraction of $E$.
Then
\[ \pi^{*} : N^1(Y_E) \longrightarrow E^{\perp} \cap N^1(Y) \]
is an isomorphism and induces an isomorphism of convex cones,
\[ \pi^{*} (\Nef(Y_E)) = \Nef(Y) \cap E^{\perp} . \]
\end{lem}

\begin{proof}
We have $N^1(Y) = \pi^{*}(N^1(Y_E)) \oplus \Z E$.
We may identify $\Nef(Y_E)$ with $\pi^{*}(\Nef(Y_E)) \subset E^{\perp}$.
The inclusion $\pi^{*}(\Nef(Y_E)) \subset \Nef(Y)$ follows immediately from
the projection formula.
This proves $\pi^{*}(\Nef(Y_E)) \subset \Nef(Y) \cap E^{\perp}$.

For the reverse inclusion, let $D \in \Nef(Y) \cap E^{\perp}$.
Since $E^{\perp} = \pi^{*}(N^1(Y_E))$, we have $D = \pi^{*} \pi_{*} D$.
Again by the projection formula, for any curve $C \subset Y_E$,
\[ \inters{\pi_* D}{C}_{Y_E} = \inters{D}{\pi^{*}C} \geq 0, \]
since $D \in \Nef(Y)$.
\end{proof}

We now prove the first of our main theorems.
We repeat it here for the convenience of the reader:
\begin{varthm}[Theorem~\ref{thm:inductive}]
Let $Y$ be a split generalized Del Pezzo surface of degree $d \leq
7$. For each $E$ in the set $\money$ of $(-1)$-curves on $Y$, let
$Y_E$ denote the split generalized Del Pezzo surface of degree $d
+ 1$ obtained by contracting $E$.  Then
\begin{equation*}
\alpha(Y) = \sum_{E \in \money} \frac{1}{d(9-d)} \alpha(Y_E).
\end{equation*}
\end{varthm}

\begin{proof}
We follow the argument used in \cite[Theorem~4]{derenthal-alpha}.
Let $\mathcal{E}$ be the set of $(-1)$- and $(-2)$-curves on $Y$.
Then $\mathcal{E}$ is exactly the set of generators for $\Eff^1(Y)$
described in Theorem~\ref{thm:effective cone generators}.
Recall that the hyperplane $\mathcal{H}_Y$ is defined as
\[
  \mathcal{H}_Y = \{ C \in N_1(Y)_{\R} : \inters{-K_Y}{C} = 1 \} .
\]
The intersection $\mathcal{P}_Y = \Nef(Y) \cap \mathcal{H}_Y$ is a polytope
with faces corresponding to $E \in \mathcal{E}$. For $E \in \mathcal{E}$, let
$\mathcal{P}_E \subset \mathcal{H}_Y$ be the convex hull of the vector
$\frac{1}{d} (-K_Y)$ and the face $\mathcal{P}_Y \cap E^{\perp}$.  (Note that
$-K_Y$ is nef by the definition of generalized Del Pezzo surface and
$\frac{1}{d}(-K_Y)$ is in $\mathcal{P}_Y$ since $\inters{-K_Y}{-K_Y} = d$.)
Then
\begin{equation*}
  \mathcal{P}_Y = \Nef(Y) \cap \mathcal{H}_Y = \bigcup_{E \in \mathcal{E}}
  \mathcal{P}_E .
\end{equation*}
The intersection of any two of the $\mathcal{P}_E$ has volume zero
in $\mathcal{H}_Y$ because the intersection lies in a subspace
of dimension strictly less than that of $\mathcal{H}_Y$.
Therefore,
\begin{equation*}
\alpha(Y) = \Vol(\mathcal{P}_Y) = \sum_{E \in \mathcal{E}} \Vol(\mathcal{P}_E).
\end{equation*}
For each $(-2)$-curve $E$, $\inters{K_Y}{E} = 0$
and hence $\frac{1}{d}(-K_Y) \in E^{\perp}$.
Thus $\mathcal{P}_E$ lies in
the hyperplane $\mathcal{H}_Y \cap E^{\perp}$
of dimension $\dim(\mathcal{H}_Y) - 1$,
and so $\mathcal{P}_E$ has volume zero.
We thus reduce to
\begin{equation*}
\Vol(\mathcal{P}_Y) = \sum_{E \in \money} \Vol(\mathcal{P}_E).
\end{equation*}

For $E \in \money$, let $\pi_E: Y \to Y_E$ be the contraction. By
Lemma~\ref{lem:nef cone identification} we have $\pi_E^* \mathcal{H}_{Y_E} =
\mathcal{H}_Y \cap E^{\perp}$. This identifies the base of the cone
$\mathcal{P}_E$ as $\mathcal{P}_Y \cap E^{\perp} = \pi_E^* \mathcal{P}_{Y_E}$.
Thus $\mathcal{P}_E$ is a cone of dimension $9-d$ with height $\frac{1}{d}$
and base volume $\Vol(\pi_E^{*} \mathcal{P}_{Y_E})$.  By Lemma~\ref{lem:nef
  cone identification}, the sublattices $N^1(Y_E) \subset N^1(Y_E)_{\R}$ and
$\pi^{*}(N^1(Y_E)) = E^{\perp} \cap N^1(Y) \subset E^{\perp}$ are isomorphic,
so $\pi^{*}$ is volume-preserving and $\Vol(\pi_E^{*} \mathcal{P}_{Y_E}) =
\Vol(\mathcal{P}_{Y_E}) = \alpha(Y_E)$.  Consequently,
\begin{equation*}
  \Vol(\mathcal{P}_E) = \frac{1}{d(9-d)} \Vol(\mathcal{P}_{Y_E}) =
  \frac{1}{d(9-d)} \alpha(Y_E) .
\end{equation*}
Summing over $E \in \money$ gives the desired result.
\end{proof}

\begin{rem}
  This generalization explains why Theorem~\ref{thm:derenthal} does not hold
  for $d=7$.  When one blows down a $(-1)$-curve on an ordinary Del Pezzo
  surface of degree $d$ for $d \leq 7$ the result is an ordinary Del Pezzo
  surface of degree $d+1$. For $d \leq 6$, the resulting ordinary Del Pezzo
  surfaces all have the same nef cone volume. This is no longer true when
  $d=7$. Let $X_d$ denote an ordinary Del Pezzo surface of degree $d$ obtained
  by blowing up $9-d$ points in general position on $\PP^2$.  Recall that $X_7
  = \Bl_{p,q}(\PP^2)$ contains three $(-1)$-curves: the exceptional divisors
  $E_p$ and $E_q$, and the proper transform $L_{pq}$ of the line through $p$
  and $q$.  Contracting $E_p$ or $E_q$ results in an $X_8$, while contracting
  $L_{pq}$ results in $\PP^1 \times \PP^1$. We have
  \[ \alpha(X_7) = \frac{1}{14} ( 2\alpha(X_8) + \alpha(\PP^1 \times \PP^1)) =
  \frac{1}{24} \] since $\alpha(X_8) = 1/6$ and
  $\alpha(\PP^1\times\PP^1)=1/4$.
\end{rem}

\section{Root systems and Weyl groups}\label{section:weyl}

In this section, we recall some of the basic facts about the root system of
$(-2)$-classes on a Del Pezzo surface and its associated Weyl group. We use
this structure in our second main result which relates the nef cone volumes of
split generalized and ordinary Del Pezzo surfaces of the same degree.

\subsection{Root systems}

\begin{defn}
A \defining{root system} $R$ is a finite collection of non-zero
vectors in a finite-dimensional real vector space $V$ with a
non-degenerate definite inner product $\inters \cdot \cdot$
satisfying the following conditions.
\begin{enumerate}
\item The set $R$ spans $V$, i.e. $R$ is essential.
\item For each $x \in R$, let $s_{\rt} : V \to V$ be the reflection through
  the hyperplane orthogonal to $x$:
  \[ s_{\rt}(v) = v - 2 \frac{\inters{\rt}{v}} { \inters{\rt}{\rt} } \rt . \]
  For each $\rt \in R$, it is required that $s_{\rt}$ takes $R$ to $R$.
\item For every $\rt_1, \rt_2 \in R$,
  \begin{equation*}
    2\frac{\inters{\rt_1}{\rt_2}}{\inters{\rt_2}{\rt_2}}
  \end{equation*}
  is an integer, i.e.,~$R$ is crystallographic.
\item If $\rt \in R$ and $c\rt \in R$, then $c \in \{ 1, -1 \}$, i.e.,~$R$ is
  reduced.
\end{enumerate}
\end{defn}

\begin{defn}
A morphism of root systems from $R \subset V$ to $R' \subset V'$
is a linear map $\Phi : V \rightarrow V'$ such that
\begin{enumerate}
\item $\Phi(R) \subset R'$, and
\item $\Phi$ preserves inner products up to a scalar multiple, i.e., there is
  a $c \in \R$ such that $\inters{\Phi(x)}{\Phi(y)} = c\cdot \inters{x}{y}$.
  Equivalently, the integers $2 \inters{x_1}{x_2}/\inters{x_2}{x_2}$ are
  preserved for all $x_1,x_2 \in R$.
\end{enumerate}
\end{defn}

\begin{rem}
We will sometimes refer to a root system $R$ in a vector space $V$
even when $R$ does not span $V$.  Strictly speaking, $R$ is only a
root system in the subspace it spans, but this minor abuse
of language should not cause any confusion.
\end{rem}

We recall some standard notions; for details,
see~\cite[Section~1.3]{humphries}, \cite[Section~VI.1.2]{bourbaki-lie-4},
\cite[Chapter~8]{MR1997306}.  Any hyperplane in $V$ not containing any root of
$R$ divides $R$ into two subsets, of \defining{positive roots} on one side
(and negative roots on the other side).  Those positive roots which cannot be
written as a sum of other positive roots with positive coefficients form a set
of \defining{simple roots}.  Each set of simple roots (for each choice of a
set of positive roots) is a linearly independent set such that every root in
$R$ is either a sum of simple roots with non-negative coefficients, or a sum
of simple roots with non-positive coefficients.

A decomposition of $R$ is a disjoint union $R = R_1 \cup \dots
\cup R_k$ such that the span of $R$ is the direct sum of the spans
of the $R_j$, each $R_j$ is a root system in its span, and the
spans of the $R_j$ are orthogonal to each other. If $R$ admits no
non-trivial decomposition, then $R$ is an \defining{irreducible
root system}. If $R$ is reducible, it has
a unique (up to order) decomposition into irreducible root systems,
called the \defining{irreducible components} of $R$.

Recall the classification of root systems by \defining{Dynkin diagrams}.  For
a root system $R$ and a choice of a set $R_0$ of simple roots in $R$, the
Dynkin diagram of $R$ is the graph with vertex set $R_0$ and an edge joining
two vertices if and only the corresponding roots are not perpendicular.  One
labels the edges of the graph according to the angle between the roots and
their relative length; for details, see~\cite{bourbaki-lie-4}.  The Dynkin
diagram is independent of the choice of a set of simple roots.  The
irreducible root systems correspond to connected graphs.  The irreducible
components of a reducible root system $R$ correspond exactly to the connected
components of the Dynkin diagram of $R$.  One has the well-known
classification of irreducible root systems corresponding to Dynkin diagrams of
types $\bA_n$, $n\geq1$; $\bB_n$, $n\geq2$; $\bC_n$, $n\geq3$; $\bD_n$,
$n\geq4$; $\bE_n$, $6\leq n \leq 8$; $\bF_4$; $\bG_2$.

The group of orthogonal transformations generated by all $s_{\rt},
\rt \in R,$ is finite and is called the \defining{Weyl group}
$W(R)$. A \defining{wall} in $V$ is a hyperplane orthogonal to an
$\rt \in R$. Removing the walls from $V$ leaves a finite set of
open convex cones called \defining{chambers}. The action of $W(R)$
permutes these chambers simply transitively.

Table \ref{table:weyl} lists all of the \defining{simply laced}
root systems (those in which all roots have the same self-intersection)
and the orders of their Weyl groups.
Table \ref{table:weyl2} gives the same data
for the non-simply laced root systems.

\begin{table}[ht]
  \centering
  \[\begin{array}{|c|c|c|c|c|c|}
    \hline
    \text{root system $R$} & \bA_n & \bD_n & \bE_6 & \bE_7 & \bE_8 \\
    \hline
     \#W(R) & (n+1)! & 2^{n-1} \cdot n! & 2^7 \cdot 3^4 \cdot 5 &
     2^{10} \cdot 3^4 \cdot 5 \cdot 7 & 2^{14} \cdot 3^5 \cdot 5^2
     \cdot 7\\
    \hline
  \end{array}\]
  \caption{The orders of simply laced Weyl groups}
  \label{table:weyl}
\end{table}

\begin{table}[ht]
  \centering
  \[
  \begin{array}{|c|c|c|c|c|}
  \hline
  \text{root system $R$} & \bB_n & \bC_n & \bF_4 &\bG_2 \\ \hline
  \#W(R) & 2^n \cdot n! & 2^n \cdot n! & 2^7 \cdot 3^2 & 2^2 \cdot 3 \\ \hline
  \end{array}
  \]
  \caption{The orders of non-simply laced Weyl groups}
  \label{table:weyl2}
\end{table}

\subsection{Root systems on Del Pezzo surfaces}\label{section:roots-on-dps}

Let $Y$ be a split generalized Del Pezzo surface of degree $d \leq 7$.  By
\cite[Sections~23--25]{manin}, the finite set $R_d$ of $(-2)$-classes on $Y$
is a root system in $N^1(Y)_{\R}$ and of course depends only on the degree
$d$. For $d \le 6$, the roots span the hyperplane $(-K_Y)^{\perp}$. The
classification of this root system is shown in Table~\ref{table:root-systems}.
\begin{table}[ht]
\[\begin{array}{| c | c | c | c | c | c | c | c | }
\hline
d & 7 & 6 & 5 & 4 & 3 & 2 & 1 \\
\hline
R_d & \bA_1 & \bA_1 \times \bA_2 & \bA_4 & \bD_5 & \bE_6 & \bE_7 & \bE_8 \\
\hline
\end{array}\]
\caption{Classification of root systems $R_d$}\label{table:root-systems}
\end{table}

Not only is $R_d$ a root system, but in fact
the subset of
$(-2)$-classes that are effective on $Y$ gives rise to a root
system, \cite[Theorem~III.2.b]{Demazure}:
\begin{thm}[Demazure]\label{thm:demazure}
  Let $Y$ be a split generalized Del Pezzo surface of degree $d \leq 6$ and
  let $\RYeff$ be the set of effective $(-2)$-classes on $Y$.  Then
  $R_Y :=\RYeff \cup -\RYeff$ is a root system in $N^1(Y)$ whose simple roots
  are the $(-2)$-curves of $Y$ and whose positive roots are $\RYeff$.
  It is contained in $R_d$.\qed
\end{thm}

\begin{rem}
  Urabe~\cite[Main Theorem]{Urabe} has shown that every root system contained
  in $R_d$ occurs as the root system $R_Y$ of a generalized Del Pezzo surface
  $Y$ of degree $d$ as in Theorem~\ref{thm:demazure}, with four exceptions:
  the subsystem of type $7 \bA_1$ in $R_2$ and the subsystems of type $7
  \bA_1$, $8 \bA_1$, and $\bD_4 + 4\bA_1$ in $R_1$.
\end{rem}

\begin{rem}\label{rem:dynkin diag}
The root system $R_Y$ can have irreducible components of the following types:
\[\bA_1, \dots, \bA_8, \, \bD_4, \dots, \bD_8, \, \bE_6, \, \bE_7, \, \bE_8.\]

For $Y$ of degree $d \geq 3$, consider the anticanonical morphism $\phi$
defined by the linear series $\lvert -K_Y \rvert$ which maps $Y$ to a
projective space of dimension $d$. For $d = 2$ (resp.~$d = 1$), let $\phi$ be
the morphism defined by the linear series $\lvert -2K_Y \rvert$ (resp.~$\lvert
-3K_Y \rvert$). Let $Y'$ be the image of $Y$ under $\phi$. The map $\phi$
sends the union of $(-2)$-curves corresponding to any connected component of
the Dynkin diagram to a singularity of $Y'$, while it is an isomorphism
between the complement of the $(-2)$-curves on $Y$ and the complement of the
singularities on $Y'$. Each singularity on $Y'$ is a rational double point.
Its type in the $\bA\bD\bE$-classification is given by the type of the
corresponding irreducible Dynkin diagram. The surface $Y'$ is a singular Del
Pezzo surface, whose minimal desingularization is the generalized Del Pezzo
surface $Y$.
\end{rem}

\subsection{Weyl groups and nef cone volume}
\label{section:weyl-groups-volume}

We proceed with the proof of our second main result,
which we repeat here for the convenience of the reader.

\begin{varthm}[Theorem~\ref{thm:weyl}]
Let $Y$ be a split generalized Del Pezzo surface of degree $d \leq
7$ and let $X$ be a split ordinary Del Pezzo surface of the same
degree. Then
\begin{equation*}
\alpha(Y) = \frac{1}{\#W(R_Y)} \alpha(X),
\end{equation*}
where $W(R_Y)$ is the Weyl group of the root system $R_Y$ whose
simple roots are the $(-2)$-curves on $Y$.
\end{varthm}

\begin{proof}
Identify $N^1(X)$ and $N^1(Y)$ as in Lemma~\ref{lem:identification}.

With notation as in the statement of Theorem~\ref{thm:demazure}, let $C$ be
the open convex cone in $N_1(Y)_{\R}$ dual to the cone spanned by the
$(-2)$-curves of $Y$.
That is, \[C = \{ \, v \in N_1(Y)_{\R} : \inters{v}{\rt} > 0 \text{ for all
  $(-2$)-curves $\rt$ on $Y$}\}.\]  Since the $(-2)$-curves are a
system of simple roots of $R_Y$,
$C$ is a single chamber for the Weyl group
$W(R_Y)$.
Recall that by Corollary~\ref{cor:cone containments},
$\Nef(Y) = \Nef(X) \cap \overline{C}$.
Intersecting
with the hyperplane $\mathcal{H}_X$ gives
$ \mathcal{P}_Y = \overline{C} \cap \mathcal{P}_X $.

We have
\[ N^1(X)_{\R} = \bigcup_{w \in W(R_Y)} \overline{ w C } , \]
so
\[ \mathcal{P}_X = \bigcup_{w \in W(R_Y)} \left( \overline{w C} \cap
  \mathcal{P}_X \right) . \] The sets $\overline{w C} \cap \mathcal{P}_X$, $w
\in W(R_Y)$, are pairwise disjoint except along boundaries, which have zero
volume.  The action of $W(R_Y)$ preserves volume and fixes $\Nef(X)$ and
$-K_X$.  Therefore it fixes $\mathcal{P}_X$, and we have
\[
\begin{split}
  \alpha(X) &= \Vol( \mathcal{P}_X )= \sum_{w \in W(R_Y)} \Vol( \overline{w C}
  \cap \mathcal{P}_X )= \#(W(R_Y)) \cdot \Vol( \overline{C} \cap
  \mathcal{P}_X ) \\
  &= \#(W(R_Y)) \cdot \Vol( \mathcal{P}_Y ) = \#(W(R_Y)) \cdot \alpha(Y).
\end{split}
\]
This completes the proof.
\end{proof}

\begin{rem}\label{rem:weyl-group-orders}
  As in Remark~\ref{rem:dynkin diag}, let $Y'$ be the singular Del Pezzo
  surface whose minimal desingularization is $Y$. The number $\#W(R_Y)$, and
  therefore $\alpha(Y)$, can be determined directly from the types of
  singularities on $Y'$ as follows: The types $R$ of the singularities of $Y'$
  coincide with the types of the irreducible components of $R_Y$. The numbers
  of elements of their Weyl groups $W(R)$ can be found in
  Table~\ref{table:weyl}. Their product is $\#W(R_Y)$.
\end{rem}

\section{Non-split Generalized Del Pezzo Surfaces}\label{section:galois}

We recall some facts about the geometry of generalized Del Pezzo
surfaces that are not split and then introduce the notion of orbit
root systems.  The results collected here will be used in Section
\ref{section:final} to relate the nef cone volume of non-split
generalized Del Pezzo surfaces to the nef cone volume of ordinary
Del Pezzo surfaces.

\subsection{The Galois action}

Throughout this section, we let $Y$ be a generalized Del Pezzo surface of
degree $d \leq 7$ defined over a field $\KK$ of characteristic $0$ and we assume
that $Y$ contains a $\KK$-rational point; we let $\overline{Y} = Y
\times_{\KK} \overline{\KK}$. The Galois group $\GalK = \Gal(\overline{\KK} /
\KK)$ acts on $N^1(\overline{Y})$, and each automorphism of
$N^1(\overline{Y})$ induced by an element of $\GalK$ preserves both the
intersection form and the anticanonical class.

\begin{prop}\label{prop:automorphisms of N^1}
The group of automorphisms of $N^1(\overline{Y})$ which preserve
the intersection form $\inters \cdot \cdot$ and the anticanonical
class $-K_Y$ is canonically isomorphic to $W(R_d)$.
\end{prop}
\begin{proof}
The result for ordinary Del Pezzo surfaces can be found in
\cite[Theorem~23.9]{manin}. (The statement there is given only for
$d\leq 6$, but the $d=7$ case is an easy calculation.) The result
holds for generalized Del Pezzo surfaces via the identification
described in Lemma~\ref{lem:identification}.
\end{proof}

Thus the action of $\GalK$ factors through (a subgroup of) the
finite group $W(R_d)$.

\begin{prop}\label{prop:galois descent on neron-severi group}
Let $Y$ be a generalized Del Pezzo surface defined over the field
$\KK$ containing a $\KK$-rational point. Then $N^1(Y) =
N^1(\overline{Y})^{\GalK}$.
\end{prop}

Recall that if $S$ is a set on which the group $G$ acts,
the standard notation $S^G = \{ s : gs=s, \forall g \in G \}$
denotes the set of fixed points of the action.

\begin{proof}
  A result of Colliot-Th{\'e}l{\`e}ne and Sansuc \cite[Theorem~2.1.2,
  Claim~(iii)]{MR899402} assures that under the hypotheses of the proposition,
  $\Pic(Y) = \Pic(\overline{Y})^{\GalK}$.  Since the intersection form on
  $\Pic(\overline{Y})$ is non-degenerate, we have $\Pic(\overline{Y}) =
  N^1(\overline{Y})$.  Finally, to show $N^1(Y) = \Pic(Y)$ it suffices to
  prove that a divisor is numerically trivial on $\Pic(\overline{Y})$ if it is
  numerically trivial on $\Pic(Y)$.

  Suppose $D \in \Pic Y$ is numerically trivial in $\Pic Y$. Let $E$ be any
  divisor class on $\overline{Y}$.  Recall that the action of $\GalK$ on
  $N^1(\overline{Y})$ factors through the finite Weyl group $W(R_d)$, so the
  $\GalK$-orbit of $E$ is finite. Say this orbit is $\{E_1, \dots, E_s\}$.
  Since $\GalK$ preserves the intersection form on $\overline{Y}$ and $D$ is
  $\GalK$-invariant,
  \[ \inters{D}{E} = \frac{1}{s} \sum_i \inters{D}{E_i} = \inters{D}{
    \frac{1}{s} \sum_i E_i} = 0 \] because $(1/s) \sum E_i$ lies in $(\Pic
  \overline{Y})^{\GalK} = \Pic Y$.

  Putting together the above results, $N^1(Y) = \Pic(Y) =
  \Pic(\overline{Y})^{\GalK} = N^1(\overline{Y})^{\GalK}$, proving the
  proposition.
\end{proof}

We now explain the relation between the effective cone of $Y$ and that of
$\overline{Y}$.

\begin{prop}\label{prop:eff cone nonsplit case}
  The effective cone of $Y$ is equal to the cone of
  $\GalK$-invariant effective classes of $\overline{Y}$, i.e.,
\begin{equation*}
  \Eff^1(Y) = \Eff^1(\overline{Y})^{\GalK}
\end{equation*}
\end{prop}

\begin{proof}
  By Proposition~\ref{prop:galois descent on neron-severi group}, we have
  $N^1(Y) = N^1(\overline{Y})^{\GalK}$.  It is clear that $\Eff^1(Y) \subseteq
  \Eff^1(\overline{Y})^{\GalK}$.  To show the reverse inclusion, first note
  that if $D$ is any effective divisor on $\Eff^1(\overline{Y})$, letting
  $\LL$ be a finite Galois extension of $\KK$ over which $D$ is defined, then
  $\sum_{\sigma \in \Gal(\LL/\KK)} \sigma(D) \in \Eff^1(Y)$.  For any $D \in
  \Eff^1(\overline{Y})^{\GalK}$ that is defined over a finite Galois extension
  $\LL / \KK$, we have
  \begin{equation*}
    D = \frac{1}{\# \Gal(\LL/\KK)} \sum_{\sigma \in \Gal(\LL/\KK)}
    \sigma(D).
  \end{equation*}
  This completes the proof.
\end{proof}

The action of $\GalK$ on $N^1(\overline{Y})$ induces an action both on the set
of $(-1)$-curves and on the set of $(-2)$-curves.

\begin{cor}\label{cor:eff cone gen nonsplit case}
  A set of generators for the effective cone of $Y$ consists of,
  for each orbit of $\GalK$ on the sets of $(-1)$-curves and $(-2)$-curves,
  the sum of the classes in that orbit.
\qed
\end{cor}

Note that this set of generators may fail to be minimal.  (See rows 3, 6 and 9
of Table \ref{table:deg5} for examples.)

\subsection{Orbit Root Systems}

We will use the following construction in Section~\ref{subsection:volume of
  pairs} in the case of $\GalK$ acting on the root system $R_{\overline{Y}}
\subset N^1(\overline{Y})$ (as in Theorem~\ref{thm:demazure}) of
$\overline{Y}$, in order to obtain a root system in $N^1(Y)$.

\begin{defn}\label{dfn:orbit root system}
Let $R \subset V$ be a possibly reducible root system with a
chosen set $\Pi$ of positive roots.  Suppose a group $G$ acts
linearly on $V$ in such a way that it permutes the elements of
$R$, preserves the inner product between elements of $R$ and
preserves positivity.  In this case, we say that $G$ acts on $R$.

The set
\begin{equation*}
\orbsys{R}{G} := \left\{ \sum_{x \in \orbit} x : \orbit \text{ is
a $G$-orbit of an element of $R$} \right\}
\end{equation*}
is called the \defining{orbit root system} of $R$ with respect to
$G$.  (We show below that $\orbsys{R}{G}$ is indeed a root
system.)
\end{defn}

\begin{prop}\label{prop:irred orbit root system}
Let $R \subset V$ be an irreducible root system with a chosen
positive system $\Pi$.
Suppose $G$ acts on $R$.
Then
$\orbsys{R}{G}$
is an irreducible root system as in Table~\ref{table:orbit root
systems}. The simple (resp. positive) roots of $\orbsys{R}{G}$ are
the sums of elements of orbits of simple (resp. positive) roots of
$R$.
\end{prop}

\begin{proof}
  Any group action which preserves inner products and positivity must
  necessarily act as an automorphism of the Dynkin diagram.  Indeed, the group
  takes non-simple roots to non-simple roots, and thus takes simple roots to
  simple roots.  Thus the group acts on the vertices of the Dynkin diagram;
  since the edges (and edge labelings) are determined by the inner product,
  they are preserved by the group.  We check case by case that all non-trivial
  admissible group actions on irreducible Dynkin diagrams are listed in
  Table~\ref{table:orbit root systems}. In each case, a direct calculation
  shows that $\orbsys{R}{G}$ is indeed a root system of the listed type.
\end{proof}

We note that a list similar to Table~\ref{table:orbit root systems} has been
compiled by Kac in \cite[Propositions~7.9, 7.10]{Kac}.  The main difference
between our list and Kac's is that we use the sum of roots in an orbit, while
he uses the average; because of this difference Kac's approach sometimes gives
the dual root system to ours.

\begin{table}
\[\begin{array}{| c | c | c |}
\hline R & G & \orbsys{R}{G} \\
\hline \bA_{2n} & \mathbb{Z} / 2\mathbb{Z} & \bB_{n} \\
\hline \bA_{2n+1} & \mathbb{Z} / 2\mathbb{Z} & \bB_{n+1} \\
\hline \bD_{n} & \mathbb{Z} / 2\mathbb{Z} & \bC_{n-1} \\
\hline \bD_{4} & \text{$\mathbb{Z} / 3\mathbb{Z}$ or
$\mathfrak{S}_3$} & \bG_2
\\
\hline \bE_{6} & \mathbb{Z} / 2\mathbb{Z} & \bF_4 \\
\hline
\end{array}\]
\caption{Non-trivial irreducible orbit root
systems}\label{table:orbit root systems}
\end{table}

\begin{lem}
Let $R \subset V$ be a possibly reducible root system with a
chosen positive system $\Pi$. Suppose $G$ acts
on $R$.
Then $G$ acts on the irreducible components of $R$ in the
following sense. Let $R = \bigcup_{i=1}^n R_i$ be a decomposition
of $R$ into irreducible components. Then for any $i$ and any $g
\in G$, $g(R_i)$ is one of the irreducible components $R_j$.
\end{lem}
\begin{proof}
One way to see this is by considering the Dynkin diagram $D$ of
$R$. Each component $R_i$ corresponds to a connected component of
the graph $D$. As noted above, the group $G$ acts as a graph
automorphism of $D$. Then each element of $G$ must take connected
components of $D$ to connected components.
\end{proof}

To avoid confusion between the actions of $G$ on $R$ and on the set of
irreducible components of $R$, we refer to orbits in the latter set as
``component orbits''.

\begin{prop}\label{prop:orbit root system}
Let $R \subset V$ be a possibly reducible root system with a
chosen positive system $\Pi$. Suppose $G$ acts
on $R$.
Let $R_1, \dots, R_k$ be irreducible components of $R$
which form a set of component orbit representatives,
that is, each component orbit contains exactly one of the $R_i$.
For each $i$,
let $G_i \subset G$ be the subgroup fixing $R_i$. Then
$\orbsys{R}{G}$ is a root system and
\begin{equation} \tag{$*$} \label{eqn:orbit root system isomorphism}
  \orbsys{R}{G} \cong \bigcup_{i=1}^k \orbsys{R_i}{G_i} .
\end{equation}
\end{prop}

\begin{proof}
First, note the right-hand side is indeed a root system. For by
Proposition~\ref{prop:irred orbit root system}, each
$\orbsys{R_i}{G_i}$ is a root system contained in the subspace
spanned by $R_i$ (since each element of $\orbsys{R_i}{G_i}$ is a
sum of one or more elements of $R_i$). Then if $i \neq j$, by
assumption $R_i$ and $R_j$ are distinct irreducible components of
$R$, so they span perpendicular subspaces of $V$. Therefore
$\orbsys{R_i}{G_i}$ and $\orbsys{R_j}{G_j}$ are perpendicular.
Hence the union on the right-hand side of~\eqref{eqn:orbit root
system isomorphism} is a perpendicular union of root systems.

Now, the spans of the component orbits are pairwise perpendicular, so we may
treat them separately.  We consider the orbit $i=1$, the others being similar.
Let the component orbit of $R_1$ consist of the components $R_{1,1}=R_1,
R_{1,2}, \dots, R_{1,p}$.  Choosing elements $g_1 = \id_G, g_2, \dots, g_p \in
G$ such that $g_i R_1 = R_{1,i}$ for each $i$, we get isomorphisms
\[ (\linspan R_1, R_1) \cong (\linspan R_{1,2}, R_{1,2}) \cong \cdots \cong
(\linspan R_{1,p}, R_{1,p}).\] Under this identification we have an
isomorphism of the diagonal
\[ \Delta \subset (\linspan R_1)^p \cong (\linspan R_{1,1}) \oplus \dots
\oplus (\linspan R_{1,p}) \]
with $\linspan(R_1)$
by projection onto the first factor.  Note that
this projection preserves angles and ratios of lengths, but divides all
lengths by a factor of $\sqrt{p}$.
One can check that the
projection takes
\[ \orbsys{R_{1,1} \cup \dots \cup R_{1,p}}{G} \]
to $\orbsys{R_1}{G_1}$, as desired.

More precisely, if $\orbit$ is the orbit of $r \in R_1$ under $G_1$, then
$\orbit\cup g_2\orbit \cup\dots \cup g_p\orbit$ is the orbit of $r$ under $G$.
Then $g_i\sum_{x \in \orbit} x = \sum_{x \in g_i\orbit} x$ is an element of
$\orbsys{R_{1,i}}{g_iG_1g_i^{-1}}$ where $g_iG_1g_i^{-1}$ is the subgroup of
$G$ fixing $R_{1,i}$, while $\sum_{i=1}^p g_i \sum_{x \in \orbit} x$ is an
element of $\orbsys{R_{1,1} \cup \dots \cup R_{1,p}}{G}$, which lies in
$\Delta$. It is projected to the element $\sum_{x \in \orbit} x$ of
$\orbsys{R_1}{G_1}$.
\end{proof}

\begin{cor} \label{cor:weyl group of orbit root system}
In the setting of Proposition~\ref{prop:orbit root system},
\[ W(\orbsys{R}{G}) \cong \prod_{i=1}^k W(\orbsys{R_i}{G_i}) . \]
\qed
\end{cor}

\section{Nef cone volume of non-split generalized Del Pezzo
surfaces}\label{section:final}

Let $Y$ be a non-split generalized Del Pezzo surface of degree at most 7,
defined over a field $\KK$ of characteristic $0$.  As in Section
\ref{section:galois} we continue to assume that Y contains a $\KK$-rational
point.  Then $\GalK = \Gal(\overline{\KK}/\KK)$ acts on the set of
$(-2)$-curves on $\overline{Y}$ and on the associated root system. In this
situation, we can construct an orbit root system as in
Definition~\ref{dfn:orbit root system}. As in the split case
(Theorem~\ref{thm:weyl}), this allows us to relate the nef cone volume of $Y$
to a volume associated to an ordinary Del Pezzo surface of the same degree.
In Section~\ref{subsection:high degree}, we compute this volume for all
non-split Del Pezzo surface of degree at least 5.

\subsection{Nef cone volume of pairs}\label{subsection:volume of pairs}

Using Proposition~\ref{prop:automorphisms of N^1}, we can
associate to
a generalized Del Pezzo surface $Y$ of degree $d \le 7$
the pair $(\overline{Y}, H_Y)$, where $H_Y
\subset W(R_d)$ is the image of $\GalK$ under the homomorphism
\begin{equation*}
\GalK \rightarrow \Aut(N^1(\overline{Y}), \langle \cdot,
\cdot \rangle, -K_Y) \cong W(R_d).
\end{equation*}
Note that $\GalK$ and therefore also $H_Y$ acts on the set of $(-2)$-curves on
$\overline{Y}$ and also on the set of its $(-1)$-curves.

\begin{rem}\label{rem:realiz}
To every generalized Del Pezzo surface $Y$ over $\KK$ there is the
associated pair $(\overline{Y},H_Y)$, as described above. The
``realization problem for pairs'' is to describe which pairs
$(\overline{Y},H)$ are obtained in this manner. That is, for which
pairs $(\overline{Y},H)$, consisting of a split generalized Del
Pezzo surface $\overline{Y}$ over $\overline{\KK}$ of degree $d$
and a subgroup $H \subset W(R_d)$ acting on the set of
$(-2)$-curves, is there a $Y$ defined over $\KK$ such that
$\overline{Y} = Y \times_{\KK} \overline{\KK}$ and $H = H_Y$ is
the image of $\GalK$ in $W(R_d)$?

Corn has shown that every pair $(\overline{X},H)$, with $\overline{X}$ a split
\emph{ordinary} Del Pezzo surface of degree $6$ and $H \subset W(R_6)$
arbitrary, is realizable in the above sense~\cite[Theorem~5.1]{corn}.

We use pairs to circumvent this realization problem. This allows us to prove
comparison theorems without having to address realization (see Corollary
\ref{cor:gen-dp-weyl}).
\end{rem}

We now define the nef cone volume $\alpha(Y,H)$ of a pair $(Y,H)$ where $Y$ is
any {\it split} generalized Del Pezzo surface of degree $d \le 7$ and
$H$ is any subgroup of $W(R_d)$ that acts on the set of
$(-2)$-curves on $Y$. It follows from Lemma \ref{lem:irred of exc
classes} that $H$ also acts on the set of $(-1)$-curves.  Note
that there is no restriction on $H$ if $Y$ is ordinary.

For such a pair $(Y,H)$, define $N^1(Y,H)$ to be $N^1(Y)^H$.
Motivated by Corollary~\ref{cor:eff cone gen nonsplit case}, we
define $\Eff^1(Y,H)$ to be the cone in $N^1(Y,H)_{\R}$
generated by the sum of the classes in each orbit of $H$ acting on the
sets of $(-1)$-curves and $(-2)$-curves of $Y$.
We naturally get a dual cone
\begin{equation*}
\Nef(Y,H) := \{ C \in N^1(Y,H)_{\R} : \langle D, C \rangle
\geq 0 \, \, \forall D \in \Eff^1(Y,H) \} .
\end{equation*}

We then have the hyperplane
\begin{equation*}
\mathcal{H}_{Y,H} := \{ C \in N^1(Y,H)_{\R} : \langle C, -K_Y
\rangle = 1\}
\end{equation*}
and the polytope
\begin{equation*}
\mathcal{P}_{Y,H} := \text{Nef}(Y,H) \cap \mathcal{H}_{Y,H}.
\end{equation*}
And so we define
\begin{equation*}
\alpha(Y,H) := \text{Vol}(\mathcal{P}_{Y,H}),
\end{equation*}
with respect to the Leray measure $d\mu$ defined in the analogous
manner to the way it was defined in Section~\ref{section:nef cone volume}.

It is immediate from Proposition~\ref{prop:galois descent on neron-severi
  group} and Corollary~\ref{cor:eff cone gen nonsplit case} that if $Y$ is any
generalized Del Pezzo surface (not necessarily split), then
\begin{equation*}
\alpha(Y) = \alpha(\overline{Y},H_Y).
\end{equation*}

\begin{lem}\label{lemma:conjugation}
Assume that $Y$ is split and let $H_1, H_2$ be two conjugate
subgroups in $W(R_d)$.  Then
\begin{equation*}
\alpha(Y,H_1) = \alpha(Y,H_2).
\end{equation*}
\end{lem}

\begin{proof}

Let $w \in W(R_d)$ be such that $H_2 = w H_1 w^{-1}$. Let
$\orbit_i$, $i \in I$, denote the orbits of the $(-1)$- and
$(-2)$-classes under $H_1$. By definition, $\Eff^1(Y,H_1)$ is
generated by the sums $\sum_{D \in \orbit_i} D$, $i \in I$. A
simple calculation shows that the orbits of these classes under
$H_2$ are given by $w \orbit_i$, $i \in I$. We have
\begin{equation*}
\alpha(Y,H_1) = \Vol\left(\left\{ C \in N^1(Y,H_1)_{\R} :
\inters{-K_{Y}}{C} = 1, \inters{C}{\sum_{D \in \orbit_i} D}
\geq 0 \,\, \forall i \in I \right\}\right).
\end{equation*}
Making use of the fact that elements of $W(R_d)$ preserve the
intersection form and anticanonical class and noting that elements
of $W(R_d)$ are orthogonal transformations and thus preserve
volumes, we compute
\begin{align*}
  \alpha(Y,H_2) & = \Vol\left(\left\{ C \in N^1(Y,H_2)_{\R} :
      \inters{-K_{Y}}{C} = 1, \inters{C}{\sum_{D \in \orbit_i} wD} \geq
      0 \,\, \forall i \in I
    \right\}\right) \\
  & = \Vol\left(\left\{ C \in N^1(Y,H_2)_{\R} : \inters{-K_{Y}}{w^{-1} C}
      = 1, \inters{w^{-1} C}{\sum_{D \in
          \orbit_i} D} \geq 0 \,\, \forall i \in I \right\}\right) \\
  & = \Vol\left( w\left\{ C \in N^1(Y,H_1)_{\R} : \inters{-K_{Y}}{C} = 1,
      \inters{C}{\sum_{D \in \orbit_i} D} \geq 0 \,\,
      \forall i \in I \right\}\right) \\
  & = \alpha(Y,H_1).
\end{align*}
This completes the proof.
\end{proof}

\begin{cor}
  Let $Y_1$ and $Y_2$ be two generalized Del Pezzo surfaces of degree $d \leq
  7$, defined over a field $\KK$ of characteristic $0$, which are geometrically
  isomorphic, i.e., $\overline{Y_1} \cong \overline{Y_2}$.  Let $H_1$ and
  $H_2$ denote the images of $\GalK$ under the respective homomorphisms $\GalK
  \rightarrow W(R_d)$.  If $H_1$ and $H_2$ are conjugate in $W(R_d)$, then
  $\alpha(Y_1) = \alpha(Y_2)$. \qed
\end{cor}

We arrive at the following analogue of Theorem~\ref{thm:weyl}.  That theorem
provided a comparison between the nef cone volumes of a split generalized Del
Pezzo surface and of a split ordinary Del Pezzo surface of the same degree.
The following theorem generalizes this to the nef cone volumes of pairs.

\begin{thm}\label{thm:alpha of nonsplit generalized DP}
  Let $Y$ be a split generalized Del Pezzo surface of degree $d \le 7$, $X$ a
  split ordinary Del Pezzo surface of the same degree, and $H$ a subgroup of
  $W(R_d)$ acting on the set of $(-2)$-curves on $Y$.  Let $R_Y$ be the root
  system whose simple roots are the $(-2)$-curves on $Y$, and let
  $\orbsys{R_Y}{H}$ be the orbit root system associated to the action of $H$
  on $R_Y$ as in Definition~\ref{dfn:orbit root system}.  Then
  \begin{equation*}
    \alpha(Y,H) = \frac{\alpha(X,H)}{\#W(\orbsys{R_Y}{H})}.
  \end{equation*}
\end{thm}

\begin{proof}
  The proof of this theorem is a generalization of the argument that proves
  Theorem~\ref{thm:weyl}. Using Lemma~\ref{lem:identification}, we identify
  $N^1(X)$ and $N^1(Y)$. This gives an identification of $N^1(X,H)$ and
  $N^1(Y,H)$. As before, $\Nef(Y,H)$ is the intersection of $\Nef(X,H)$ with
  the closure of a chamber defined by the simple roots of $\orbsys{R_Y}{H}$.
  As in the proof of Theorem~\ref{thm:weyl}, the chambers of the Weyl group
  $W(\orbsys{R_Y}{H})$ intersect only along boundaries, which have zero
  volume.  They fill $N^1(Y)$.  There are $\#W(\orbsys{R_Y}{H})$ of the
  chambers.  From here, the proof is completed by the same steps as in the
  proof of Theorem~\ref{thm:weyl}.
\end{proof}

We arrive at our third main result, the computation of the nef cone volume of
a generalized Del Pezzo surface over a non-closed field of characteristic~$0$:

\begin{cor}\label{cor:gen-dp-weyl}
  Let $Y$ be a generalized Del Pezzo surface of degree $d \le 7$ over the
  field $\KK$ of characteristic $0$ and $X$ a split ordinary Del Pezzo surface
  of the same degree.  Let $\overline{Y} = Y \times_{\KK} \overline{\KK}$, and
  identify $N^1(\overline{Y})$ with $N^1(X)$ as in
  Lemma~\ref{lem:identification}.  Let $H_Y \subset W(R_d)$ be the image of
  $\GalK$.  Let $R_{\overline Y} \subset R_d$ be the root system whose simple
  roots are $(-2)$-curves on $\overline{Y}$.  Then
  \[ \alpha(Y) = \alpha(\overline{Y},H_Y) = \frac{ \alpha(X,H_Y)}{ \#W(
    \orbsys{R_{\overline{Y}}}{H_Y} )} . \] \qed
\end{cor}

Using Proposition~\ref{prop:irred orbit root system} and
Corollary~\ref{cor:weyl group of orbit root system}, the integer appearing in
the denominator is straightforward to compute.  This reduces the computation
of the nef cone volume of an arbitrary generalized Del Pezzo surface over a
non-closed field to the computation of the nef cone volume of a pair involving
a split ordinary Del Pezzo surface.

\subsection{Pairs involving ordinary Del Pezzo surfaces of high degree}
\label{subsection:high degree}

As examples, let us compute $\alpha(X)$ for the various possible
non-split ordinary Del Pezzo surfaces $X$ of degree $d \geq 5$.

For $d \geq 7$ there are very few possible non-trivial Galois actions on
$\overline{X}$, and we list these cases briefly:
\begin{enumerate}
\item In fact there are no non-trivial possibilities with $d=9$:
we must have $\overline{X} \cong \PP^2$,
the Galois action is trivial, and $\alpha(X) = 1/3$.
\item For $d=8$, the only possible non-trivial form occurs
when $X$ is a twist of $\PP^1 \times \PP^1$ in which the
Galois action permutes the two generating rulings.
In this case $\alpha(X) = 1/2$.
\item For $d=7$, the only possible non-trivial form occurs
when $X$ is the blow-up of two conjugate rational points on $\PP^2$,
so that the Galois action interchanges the points.
In this case $\alpha(X) = 1/6$.
\end{enumerate}

For $d=5,6$ there are many more cases.  For the remainder of this section, let
$X$ be a possibly non-split ordinary Del Pezzo surface of degree $5$ or $6$
defined over a non-closed field $\KK$ of characteristic $0$.  Let $\overline{X}
= X \times_{\KK} \overline{\KK}$.  As above, we have
\[ \alpha(X) = \alpha(\overline{X} , H_X) \]
where $H_X$ is the image of the Galois group in $W(R_d)$.
We compute $\alpha(X)$ by finding the values of $\alpha(\overline{X},H)$
for all subgroups $H$ of $W(R_d)$.
(As noted in Remark~\ref{rem:realiz}, it is not obvious
which $H \subset W(R_d)$ arise as images of Galois groups,
so a priori some values $\alpha(\overline{X},H)$ might not
correspond to any $\alpha(X)$.)

\begin{figure}\label{fig:fig1}
\includegraphics[height=2in]{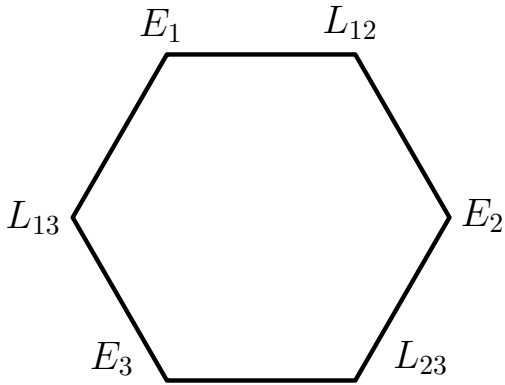}
\caption{Configuration of $(-1)$-curves on an ordinary Del Pezzo surface of degree $6$}
\label{figure:degree6}
\end{figure}

For the case $d=6$, recall that $\overline{X}$ is obtained by blowing up three
non-collinear points in $\PP^2$ and the cone $\Eff^1(\overline{X})$ is
minimally generated by the $(-1)$-curves on $\overline{X}$.  Let $E_1$, $E_2$,
$E_3$ denote the exceptional curves and $L$ denote the pullback of a line.
The set of $(-1)$-curves is shown schematically in
Figure~\ref{figure:degree6}.
In this graph, the vertices correspond to the generating classes
for $\Eff^1(\overline{X})$, with the convenient shorthand $L_{ij}
= L - E_i - E_j$.  Two classes intersect if and only if the
corresponding vertices in the graph are connected by an edge.

In Table~\ref{table:deg6}, we consider the various subgroups of
\begin{equation*}
W(R_6) = W(\bA_1) \times W(\bA_2) \cong \mathbb{Z} / 2\mathbb{Z}
\times \mathfrak{S}_3 \cong D_6.
\end{equation*}
By Lemma~\ref{lemma:conjugation}, it suffices to consider subgroups up to
conjugacy.  For each conjugacy class, we choose a representative subgroup $H$
and give the order $\#H$ of $H$, the orbit structure of $H$ on the generators
of $\Eff^1(\overline{X})$, the rank $\rho$ of $N^1(\overline{X},H)$, the
number $m$ of generators in the minimal generating set of
$\text{Eff}^1(\overline{X},H)$, and finally the nef cone volume
$\alpha(\overline{X},H)$.  We describe $H$ in terms of generators, using the
generator $s_{123}:=s_{L-E_1-E_2-E_3}$ ($180^{\circ}$ rotation) of $W(\bA_1)$
and the generators $s_{12}:=s_{E_1-E_2}$ (the flip swapping $E_1$ and $E_2$)
and $s_{23}:=s_{E_2-E_3}$ (the flip swapping $E_2$ and $E_3$) of $W(\bA_2)$.

Given $H$, we may compute $\alpha(\overline{X},H)$ as follows.  We explicitly
compute the sums of elements in each orbit of the action of $H$ on the
generators of $\Eff^1(\overline{X})$, obtaining a set of generators of the
cone $\Eff^1(\overline{X},H)$. We compute the dual cone in
$N^1(\overline{X},H)$, obtaining $\Nef(\overline{X},H)$.  Intersecting with
the hyperplane $\mathcal{H}_{\overline{X},H}$ gives the polytope
$\mathcal{P}_{\overline{X},H}$, whose volume is $\alpha(\overline{X},H)$.

\begin{table}
\[\begin{array}{| c | c | c | c | c | c |}
\hline H & \#H & \text{Orbit structure} & \rho & m & \alpha(\overline{X},H) \\
\hline \langle s_{123}, s_{12}, s_{23} \rangle & 12 & & & & \\
\langle s_{23} s_{123}, s_{23} s_{12} \rangle & 6 & & & & \\
\langle s_{123} s_{12} s_{23} \rangle & 6 & \{ E_1, E_2, E_3, L_{12}, L_{13}, L_{23} \} & 1 & 1 & 1 \\
\hline \langle s_{12}, s_{23} \rangle & 6 & & & & \\
\langle s_{12} s_{23} \rangle & 3 & \{ E_1, E_2, E_3 \}, \{ L_{12}, L_{13}, L_{23} \} & 2 & 2 & 1/3 \\
\hline \langle s_{123}, s_{23} \rangle & 4 & \{ E_1, L_{23} \},
\{E_2, E_3, L_{12}, L_{13} \} & 2 & 2 & 1/2 \\
\hline \langle s_{123} s_{12} \rangle & 2 & \{ E_1, L_{12} \}, \{ E_2, L_{13} \}, \{ E_3, L_{23} \} & 2 & 2 & 1/2 \\
\hline \langle s_{123} \rangle & 2 & \{ E_1, L_{23} \}, \{ E_2, L_{13} \}, \{ E_3, L_{12} \} & 3 & 3 & 1/8 \\
\hline \langle s_{12} \rangle & 2 & \{ E_1, E_2 \}, \{ E_3 \}, \{ L_{12} \}, \{ L_{13}, L_{23} \} & 3 & 4 & 1/12 \\
\hline \langle e \rangle & 1 & \{E_1\}, \{E_2\}, \{E_3\}, \{L_{12}\}, \{L_{13}\}, \{L_{23}\} & 4 & 6 & 1/72\\
\hline
\end{array}\]
\caption{Values of $\alpha(\overline{X},H)$ for $\overline{X}$ a split ordinary
Del Pezzo surface of degree $6$}
\label{table:deg6}
\end{table}

For the case when $d=5$, recall that $\overline{X}$ is the blowup of $\PP^2$
at $4$ points in general position.  Similarly to the case $d=6$, the cone
$\Eff^1(\overline{X})$ is generated by the $(-1)$-curves $E_i$ for $1 \leq i
\leq 4$ and $L_{ij} = L - E_i - E_j$ for $1 \leq i < j \leq 4$.

\begin{figure}
\includegraphics[height=2in]{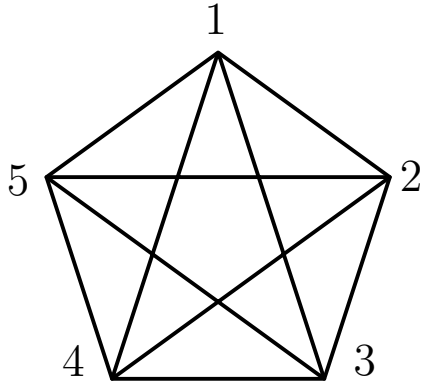}
\caption{Configuration of $(-1)$-curves on an ordinary Del Pezzo surface of degree $5$}
\label{figure:degree5}
\end{figure}

In Figure~\ref{figure:degree5}, we use a different diagram to exhibit the full
symmetry of the configuration of these $10$ curves with respect to $W(R_5) =
\mathfrak{S}_5$. (It seems impossible to make visible all of the symmetries
when drawing a diagram analogous to Figure~\ref{fig:fig1}).
Here the minimal generators of $\Eff^1(\overline{X})$ correspond to {\it
edges} of the graph, and two generating classes intersect if and only if
the corresponding edges
do {\it not} share a common vertex.  The action of $W(R_5)
= \mathfrak{S}_5$ corresponds to permuting the 5 vertices.
Table~\ref{table:gens-R5} shows the correspondence between the edges
of the diagram and the generating classes,
where we use the notation [$ij$] to indicate
the edge connecting vertex $i$ with vertex $j$.

We note that the enumeration of the conjugacy classes of subgroups of
$\mathfrak{S}_5$ has been made by G{\"o}tz Pfeiffer and is available online
\cite{Pfeiffer}.  Table~\ref{table:deg5} contains the values of
$\alpha(\overline{X},H)$ for the various possible conjugacy classes of
subgroups of $\mathfrak{S}_5$.

\begin{table}
\[\begin{array}{| c | c | c | c | c | c | c | c | c | c |}
  \hline [12] & [13] & [14] & [15] & [23] & [24] & [25] & [34] & [35] & [45] \\
  \hline E_1 & E_2 & E_3 & E_4 & L_{34} & L_{24} & L_{23} & L_{14} & L_{13}
  & L_{12}  \\
  \hline
\end{array}\]
\caption{Correspondence of edges in Figure~\ref{figure:degree5} to
generators of the effective cone of an ordinary Del Pezzo surface
of degree $5$}\label{table:gens-R5}
\end{table}

\begin{table}
\[\begin{array}{| c | c | c | c | c | c |}
\hline H & \#H & \text{Orbit structure} & \rho & m & \alpha(\overline{X},H) \\

\hline \langle (1 2), (1 2 3 4 5) \rangle & 120 & & & & \\
\langle (1 2)(3 4), (2 5 3) \rangle & 60 & & & & \\
\langle (1 2 3 4), (1 3)(2 4), (1 2 5 4 3) \rangle & 20 & \{\eu, \ed, \et, \eq, \lud, & & & \\
\langle (1 2)(3 4), (1 3 5 4 2) \rangle & 10 & \lut, \luq, \ldt, \ldq, \ltq\} & 1 & 1 & 1 \\

\hline \langle (1 2), (1 2 3 4) \rangle & 24 & \{\eu, \ed, \et, \luq, \ldq, \ltq\}, & & &\\
\langle (1 2 3), (1 2)(3 4), (1 4)(2 3) \rangle & 12 & \{\eq, \lud, \lut, \ldt\} & 2 & 2 & 2/3 \\

\hline \langle (1 2), (3 4), (3 4 5) \rangle & 12 & & & & \\
\langle (1 2)(3 4), (3 4 5) \rangle & 6 & \{\eu\}, \{\lud, \lut, \luq\} & & & \\
\langle (1 2), (3 4 5) \rangle & 6 & \{\ed, \et, \eq, \ldt, \ldq, \ltq\}, & 2 & 2 & 1/2 \\

\hline \langle (1 2), (1 2 3) \rangle & 6 & \{\eu, \ed, \ltq\}, \{\et, \luq, \ldq\}, & & & \\
\langle (1 2 3) \rangle & 3 & \{\eq, \lut, \ldt\}, \{\lud\} & 3 & 4 & 5/24 \\

\hline & & \{\eu, \eq, \lud, \luq, \ltq\}, & & & \\
\langle (1 2 3 4 5) \rangle & 5 & \{\ed, \et, \lut, \ldt, \ldq\} & 1 & 1 & 1 \\

\hline & & \{\eu\}, \{\ed, \et, \ldq, \ltq\}, & & & \\
\langle (1 2)(3 4), (1 3)(2 4) \rangle & 4 & \{\eq, \ldt\}, \{\lud, \lut\}, \{\luq\} & 3 & 4 & 1/6 \\

\hline & & \{\eu, \et, \luq, \ltq\}, \{\ed, \ldq\}, & & & \\
\langle (1 2), (3 4) \rangle & 4 & \{\eq, \lud, \lut, \ldt\} & 2 & 2 & 2/3 \\

\hline \langle (1 2), (3 4), (1 3)(2 4) \rangle & 8 & \{\eu, \luq\}, \{\ed, \ldq\}, \{\et, \ltq\}, & & & \\
\langle (1 2 3 4) \rangle & 4 & \{\eq, \lud, \lut, \ldt\} & 2 & 2 & 2/3 \\

\hline & & \{\eu\}, \{\ed,\ldq\}, \{\et,\ltq\}, & & & \\
\langle (1 2)(3 4) \rangle & 2 & \{\eq\,\ldt\}, \{\lud, \lut\}, \{\luq\} & 3 & 4 & 1/6 \\

\hline & & \{\eu\}, \{\ed,\ltq\}, \{\et,\ldq\}, & & & \\
\langle (1 2) \rangle & 2 & \{\eq\,\ldt\}, \{\lud\}, \{\lut\}, \{\luq\} & 4 & 7 & 1/24 \\

\hline & & \{\eu\}, \{\ed\}, \{\et\}, \{\eq\}, \{\lud\}, & & & \\
 \langle e \rangle & 1 & \{\lut\}, \{\luq\}, \{\ldt\}, \{\ldq\},
\{\ltq\}
 & 6 & 10 & 1/144\\
\hline
\end{array}\]
\caption{Values of $\alpha(\overline{X},H)$ for $\overline{X}$ a split ordinary
Del Pezzo surface of degree $5$}
\label{table:deg5}
\end{table}

\bibliography{dp-bib}

\end{document}